\documentclass[11pt]{article}

\usepackage{color}
\usepackage{graphicx}
 
\usepackage{url}
\usepackage{enumitem} 
\usepackage{graphicx}
\usepackage{color}
\usepackage{rotating}

\usepackage{algorithm}
\usepackage{algorithmic}

\usepackage{amssymb}
\usepackage{amsmath}
\usepackage{amsfonts}
\usepackage{amsthm}
 
 \usepackage{mathtools}
  
\usepackage{float}
\usepackage{ulem}

\usepackage{array,arydshln}

\newcommand\by {\mathbf y}

\newcommand\indica {\mathbb{I}}

\newcommand\itA {{\mathcal{A}}}
\newcommand\itB {{\mathcal{B}}}

\newcommand\itE {{\mathcal{E}}}
\newcommand\itF {{\mathcal{F}}}

\newcommand\itH {{\mathcal{H}}}
\newcommand\itI {{\mathcal{I}}}

\newcommand\itM {{\mathcal{M}}}
\newcommand\itN {{\mathcal{N}}}

\newcommand\itT {{\mathcal{T}}}
\newcommand\itU {{\mathcal{U}}}

\newcommand\bxi {\mbox{\boldmath $\xi$}}

\newcommand\bGa {\mbox{\boldmath $\Gamma$}}

\newcommand\wphi {\widehat{\phi}}

\newcommand\wlam {\widehat{\lambda}}

\newcommand\wmu {\widehat{\mu}}

\newcommand\wtheta {\widehat{\theta}}

\newcommand\wtau {\widehat{\tau}}

\newcommand\wGamma {\widehat{\Gamma}}
\newcommand\wbGa  {\widehat{\bGa}}

\newcommand\wUps  {\widehat{\Upsilon}}

\newcommand\wPi {\widehat{\Pi}}

\def\real{\mathbb{R}}
\def\natu{\mathbb{N}}

\newcommand{\esp}{\mathbb{E}}
\newcommand{\prob}{\mathbb{P}}
\newcommand{\cov}{\mbox{\sc Cov}}

\newcommand{\convpp}{ \buildrel{a.s.}\over\longrightarrow}

\newcommand{\convprob  }{ \buildrel{p}\over\longrightarrow}

\newcommand{\convdist}{ \buildrel{D}\over\longrightarrow}

\newcommand{\trasp}{^{\mbox{\footnotesize \sc t}}}

\newcommand\bcero {{\bf{0}}}

\def\argmin{\mathop{\mbox{argmin}}}

\newcommand{\identidad}{\mbox{\bf I}}

\newcommand{\boot}{\mbox{\footnotesize\sc b}}
\newcommand{\weigh}{\mbox{\footnotesize\sc w}}

\newcommand{\gaus}{\mbox{\footnotesize\sc g}}

\newcommand{\pool}{\mbox{\footnotesize\sc pool}}

\newcommand\new{\newline}
\newcommand\noi{\noindent}

\def\square{\ifmmode\sqr\else{$\sqr$}\fi}
\def\sqr{\vcenter{
         \hrule height.1mm
         \hbox{\vrule width.1mm height2.2mm\kern2.18mm
\vrule width.1mm}
         \hrule height.1mm}}

\newcommand\Signo {\mbox{\footnotesize \sc s}}

\begin{document}

\title{The spatial sign covariance operator: Asymptotic results and applications}
\author{{ Graciela Boente},  {Daniela Rodriguez} {and} { Mariela Sued}\\
{\small \sl Facultad de Ciencias Exactas y Naturales, Universidad de Buenos Aires and
CONICET, Argentina} \\
{\small e--mail: gboente@dm.uba.ar\hskip1truecm drodrig@dm.uba.ar\hskip1truecm msued@dm.uba.ar}
}
\date{}
\maketitle

\begin{abstract}
Due to the increasing recording capability, functional data analysis has become an important research topic. For functional data the study of outlier detection and/or the development of robust statistical procedures has started recently. One robust alternative to the sample covariance operator is the  sample spatial sign covariance  operator. In this paper, we study the asymptotic behaviour of the sample spatial sign covariance  operator  when
location is unknown. Among other possible applications of the obtained results, we derive the asymptotic distribution of the principal directions obtained from the  sample spatial sign covariance  operator and   we develop   test to detect  differences between the scatter  operators of two populations.
In particular, the test performance  is illustrated through a Monte Carlo study for small sample sizes.
\end{abstract}

\normalsize

\section{Introduction}

Functional data analysis is a field which deals with a sample of curves registered on a continuous period of time. A more general and inclusive framework that can accommodate the situation in which the observations are  images or surfaces is to consider realizations of a random element $X$ on a Hilbert space $\itH$ with inner product $\langle \cdot, \cdot\rangle$ and norm $\|\cdot\|$.
The area has attracted much interest in the statistics community and has increase its development, since technological advances in data collection and
storage require procedures specifically designed for dealing with such data. It has been extensively discussed that simplifying the functional model by discretizing the observations as sequences of numbers can often fail to capture some of its important characteristics, such as the smoothness and continuity of the underlying functions. Statistical  methods to analyse such functional data may be found, for instance, in Ramsay and Silverman (2005),  Ferraty and Vieu (2006), Ferraty and Romain (2010), Horv\'ath and Kokoszka (2012) and  Hsing and Eubank (2015). For a summary of recent advances in functional statistics see Cuevas (2014) and Goia and Vieu (2016). 

In this setting, the analysis of the covariance operator arises in many applied contexts. In particular, functional principal component analysis is a common tool to explore the characteristics of the data within a reduced dimensional space. As it is well known, the principal directions may be obtained as the eigenfunctions of the covariance operator.  By exploring this lower dimensional principal components space,  functional principal components analysis  allows to detect  atypical observations or outliers in the data set, when combined with a robust estimation procedure.  Among other procedures to robustly estimate the principal directions, we can mention the   spherical principal components  of Locantore \textsl{et al.} (1999) and Gervini (2008)   that correspond to the eigenfunctions of the spatial sign operator, the projection--pursuit given in Bali \textsl{et al.} (2011), the robust approach  given on    Sawant \textsl{et al.} (2012), the $M-$type smoothing spline estimators proposed in Lee \textsl{et al.} (2013) and the $S-$estimators of Boente and Salibi\'an--Barrera (2015). One key point when deriving detection rules is that the robust functional principal direction estimators are indeed estimating the target directions. In this sense, the projection--pursuit given in Bali \textsl{et al.} (2011) and the spherical   principal components  are Fisher consistent for elliptically distributed random elements. The result for the spherical   principal components, derived in Boente \textsl{et al.} (2014), extends a previous one obtained in Gervini (2008) for random elements with a finite Karhunen-Lo\`eve expansion.  Moreover,    Boente \textsl{et al.} (2014) proved that  the linear space spanned by the  $q$   eigenfunctions associated to the $q$ larger eigenvalues  of the spatial sign covariance  operator, provides the best $q$ dimensional approximation to the centered process, in the sense of having stochastically smallest residuals squared norms among all linear spaces of dimension $q$. This result does not required second order moment for the process.  Up to our knowledge, the asymptotic distribution of the robust principal direction estimators mentioned above is unknown.  One of the goals of this paper is to derive the asymptotic distribution of the spherical principal component estimators through that of the sample spatial sign covariance  operator.

A more recent statistical problem is that of testing for equality or proportionality between the covariance operators of two  populations.
For instance,  Ferraty \textsl{et al.} (2007) considered tests for comparing groups of curves based on comparison of their covariances.
By the Karhunen--Lo\'eve expansion,
this is equivalent to testing if all the samples have the same set of functional principal components sharing also their size.
When considering only two populations, Benko \textsl{et al.} (2009), Panaretos \textsl{et al.} (2010) and Fremdt \textsl{et al.} (2013)
used this characterization to develop   test statistics.
In particular, Benko \textsl{et al.} (2009) proposed two--sample bootstrap tests for specific aspects of the spectrum of functional data,
such as the equality of a subset of eigenfunctions. On the other hand, Panaretos \textsl{et al.} (2010) and Fremdt \textsl{et al.} (2013)
considered an approach based on the projection of the data over a suitable  chosen finite--dimensional space, such as that defined by the
functional principal components of each population.  The results in  Fremdt \textsl{et al.} (2013) generalized those provided in
Panaretos \textsl{et al.}  (2010) which assume that  the processes  have a Gaussian distribution. More recently, Pigoli \textsl{et al.} (2014)
developed a two--sample test for comparing covariance operators using different   distances between covariance operators. Their procedure is based on a permutation
test and assumes that the two samples have the same mean, otherwise, an approximate permutation test is considered  after the processes are centered using their
sample means. Some authors have also consider robust proposals for this  problem.  Kraus  and Panaretos  (2012) introduced
a class of dispersion operators and proposed a procedure for testing for equality of dispersion operators among two populations. 
Recently,   Boente  \textsl{et al.}  (2017),  extended the classical two populations problem,
presenting a test for equality of covariance operators among $k$ populations in which the asymptotic distribution of the sample covariance operator plays a crucial role in deriving the asymptotic distribution of the  proposed statistic.
It is well known that  the presence of  outliers in  the sample might  drive to invalid conclusions. This motivate  the development of robust procedures to deal with these kind of problems. 
In this paper,  we also present as application of our results  a test for equality of  the spatial sign covariance  operators between two populations. The statistic mimics the one presented for the classical setting  and, as in the classical setting, its asymptotic distributions depends on that of
the  empirical spatial sign  covariance operator for each population.  It is worth noticing that, for functional elliptical distributions, equality of spatial sign  covariance operators guarantees  that  the considered populations have the same principal components.

Unlike  the classical case, where the estimation of the mean plays no role in the asymptotic distribution of the covariance operator estimator, the imputation of an estimated location  when defining the spatial sign covariance estimator requires  some special considerations. One of the goals of this paper is to present  a detailed proof  of the asymptotic distribution of the sample spatial sign covariance  estimator, which extends to the functional setting the results given   by D\"urre \textsl{et al.} (2014) in the finite--dimensional case.

The paper is organized as follows. In Section \ref{cov}, we  introduce the notation to be used in the paper as well as the spatial sign covariance operator  with unknown location, while Section \ref{asym} deals with its consistency and 
asymptotic normality. Section \ref{aplication} considers the application of the obtained results to two situations: the asymptotic behaviour of  the spherical principal component estimators  and the proposal of a      test to detect  differences between the     spatial sign covariance operators of two populations,  whose performance is also numerically studied for small samples. 
Proofs are relegated to the Appendix.

\section{The spatial covariance operator}\label{cov}
Let $\itH$ be a separable Hilbert space, such as $L^2({\itI})$  for some bounded interval $\itI$,
 with inner
product $\langle \cdot,\cdot\rangle$ and norm $\|u\|=\langle u,u\rangle^{1/2}$. The functional sign of
$u \in\itH$,  is defined as $s(u)=u/\| u\|$, for  $u\not=0$, and $s(0)=0$.

Let $X$ be a random element taking values in   $\itH$.
For a given  $t\in \itH$, the spatial or sign covariance  operator of $X$ centered at  $t$
is defined by
\begin{equation}
\Gamma^{\Signo}(t)=\esp[s(X-t) \otimes s(X-t)],
\end{equation}
where $\otimes$ denotes  the tensor product on $\itH$, e.g., for $u,
v\in {\itH}$, the operator $u\otimes v:{\itH}\to \itH$ is
defined as $(u\otimes v)w= \langle v,w\rangle u$.
Notice that  $u\otimes v$ is a compact operator
that belongs to $\itF$, the Hilbert space of  Hilbert--Schmidt operators over $\itH$.
Recall that for $\Upsilon\in \itF$,  $\Upsilon^{*}$  denotes the adjoint of the operator $\Upsilon$, while for
$\Upsilon_1,\Upsilon_2\in \itF$,  the inner product in $\itF$ is defined as  $\langle \Upsilon_1,
\Upsilon_2\rangle_{\itF} =\mbox{trace} (\Upsilon_1^{*}\Upsilon_2)
=\sum_{\ell=1}^\infty \langle \Upsilon_1 u_\ell, \Upsilon_2 u_\ell\rangle$, and so  the
norm equals $\|\Upsilon\|_{\itF}=\langle \Upsilon^{*} , \Upsilon \rangle_{\itF}^{1/2}
=\{\sum_{\ell=1}^\infty \|\Upsilon u_\ell\|^2\}^{1/2}$, with $\{u_\ell :
\ell\,\ge 1\}$ any orthonormal basis of $\itH$.
These definitions are independent of the basis choice.

Given independent  random elements $X_1,\dots,X_n$, distributed as $X$,  for each $t\in \itH$  define the sample version of $\Gamma^{\Signo}(t)$ as
$$\wGamma_n^{\Signo}(t)=\frac{1}{n}\sum_{i=1}^n s(X_i-t) \otimes s(X_i-t)\,.$$
The law of large numbers in $\itF$, entails that, for any $t\in \itH$, $\wGamma_n^{\Signo}(t)$ converges almost surely to $\Gamma^{\Signo}(t)$.
Moreover, the asymptotic distribution
can be obtained from the central limit theorem in $\itF$, see, for instance, Dauxois \textsl{et al.} (1982). 

Typically, the spatial operator is centered using   as location the functional median $\mu$  of the process $X$, that is, the object of interest is the spatial operator $\Gamma^{\Signo}(\mu)$. However, in most situations   $\mu$ is unknown.  Hence, when  estimating the spatial sign operator, an estimator of $\mu$ must be considered. More precisely, let $\wmu_n$ be  a preliminary consistent estimator of $\mu$, then   $\wGamma_n^{\Signo}(\wmu_n)$ provides an estimator of $\Gamma^{\Signo}(\mu)$. The asymptotic properties of $\wGamma_n^{\Signo}(\wmu_n)$ are presented in Section \ref{asym}.


\subsection{Some general comments}{\label{comments}} 
As mentioned in the Introduction, the sample spatial operator  $\wGamma_n^{\Signo}(\wmu_n)$ has been used as an alternative to the sample covariance operator when considering  robust estimation procedures. In particular, it has been considered when one suspects that the underlying distribution may not have finite moments.  Elliptical random elements have been introduced in Bali and Boente (2009) and further studied in Boente \textsl{et al.} (2014).  For completeness, we recall their definition.  Given a random element $X$ in a separable Hilbert space $\itH$, we say that $X$ has an elliptical distribution with parameters $\mu\in \itH$ and $\Gamma:\itH \rightarrow \itH$, where $\Gamma$ is a self--adjoint, positive semi--definite and compact operator, if and only if for any linear and bounded operator $A:\itH \rightarrow \real^d$  we have that  the vector $A\,  X$ has a $d-$variate elliptical distribution with location parameter $A \, \mu $, shape matrix  $A \,  \Gamma \, A^*$ and characteristic generator $\varphi$, that is, $A\,  X\sim \itE_d(A \, \mu, A \,  \Gamma \, A^*,\varphi)$ where $A^*:\real^d\rightarrow \itH$ denotes the adjoint operator of $A$.  We write $X \, \sim \,\itE(\mu, \Gamma,\varphi)$ and $\Gamma$ is called the scatter operator of $X$. Hence, elliptical families provide a more general setting than considering  Gaussian random elements  and   the sign operator gives a useful tool to obtain Fisher--consistent estimators of the principal directions, that is,  estimators consistent to the eigenfunctions of the scatter operator of the elliptical process, even when second moments do not exist (see Boente \textsl{et al.}, 2014). 

 For elliptical random elements, two situations may arise, either the scatter operator $\Gamma$ has a finite rank $q$ or it does not have a finite rank. In the first case, the process $X$ has a finite Karhunen--Lo\`eve expansion $X=\mu + \sum_{j=1}^q \lambda_j^{1/2} \xi_j \phi_j$, where $\phi_j$ are the eigenfunctions of $\Gamma$ related to the eigenvalues $\lambda_j$ ordered in decreasing order and $\bxi=(\xi_1,\dots, \xi_q)\trasp\sim \itE_q(0, \identidad_q,\varphi) $, that is, $\bxi$ has an spherical distribution. In this setting, the asymptotic behaviour of $\wGamma_n^{\Signo}(\wmu_n)$ may be derived from the results given in D\"urre \textsl{et al.} (2014), since the distribution of $\mbox{diag}(\lambda_1^{1/2},\dots,\lambda_q^{1/2}) \bxi$ is symmetric around $\bcero$. On the other hand, if $\Gamma$ has not a finite rank, Proposition 2.1 in Boente \textsl{et al.} (2014) states that the process is a scale mixture of Gaussian distributions, more precisely there exists a zero mean Gaussian random element $Y$  and a random variable $V>0$ independent of $Y$ such that $X=\mu\,+\, V\,Y$. Without loss of generality, throughout the paper, we will assume   that $\Gamma$ is  the covariance operator of $Y$. The results given in Section \ref{asym}  include this case but they also provide a consistency and asymptotic normality results in a framework more general than elliptical families.

\section{Asymptotic results}\label{asym}

The following results establish the consistency and the asymptotic normality of the spatial sign covariance operator with unknown location.
The proofs are relegated to the Appendix. From now on,   the notation $u_n\convpp u$ in $\itH$  means that $\|u_n-u\|\convpp 0$, while for random operators $\Upsilon_n\in\itF$, the convergence $\Upsilon_n  \convpp \Upsilon$ in $\itF$ stands for $\|\Upsilon_n -\Upsilon\|_{\itF}\convpp 0$.

\noi \textbf{Theorem \ref{asym}.1.}  \textsl{Let $\mu\in \itH$ be the location parameter of the process. Assume that  $\wmu_n$ is  strongly consistent  estimator of $\mu$ and that $\esp\left[\| X-\mu\|^{-1}\right]<\infty$.  Then , we have that $ \wGamma^{\Signo}(\wmu_n)\convpp \Gamma^{\Signo}(\mu)$.
}

\vskip0.1in
\noi \textbf{Remark \ref{asym}.1.}
In  a robust context, several  estimators of the location parameter $\mu$  have been considered. Among others we can mention  the trimmed means proposed by Fraiman and  Mu\~niz (2001), the depth--based estimators of Cuevas \textsl{et al.} (2007) and L\'opez--Pintado and  Romo (2007), or the functional median defined in Gervini (2008). In particular,  as mentioned above, the spatial median is the usual choice to center the data when location is unknown and the spatial covariance operator is considered.  The spatial median is defined as 
\begin{equation}
\mu=\argmin_{u\in \itH}  \esp[\|X-u\|-\|X\|]\,.
\label{mediana}
 \end{equation}
and different methods have been proposed to provide estimators, in the functional case. Gervini (2008) shows that the sample spatial median, denoted $\wmu_n$ and defined as the solution of the empirical version  of (\ref{mediana}), can be found solving  a convex
$n-$dimensional minimization problem. Furthermore, $\wmu_n$ is strongly consistent with respect to the weak topology in $\itH$, that is for any $u\in \itH$, $\langle \wmu_n, u \rangle \convpp \langle \mu, u \rangle$.
On the other hand, Cardot \textsl{et al.} (2013) propose  to estimate the spatial median  through an algorithm  that can be seen
as a stochastic gradient algorithm.   Theorem 3.1 in Cardot \textsl{et al.} (2013)  shows that  this estimator  converges to the median almost surely, under mild conditions. This result guarantees the existence of strong consistent estimators of the median in the functional case and hence, that of the estimators of the spatial sign covariance operator.

\vskip0.1in
In order to study the asymptotic distribution of $\wGamma^{\Signo}(\wmu_n)$, let
$\itB$ denotes the Banach space of linear and continuous operators from $\itH$ to $\itF$, that is,
 $\itB=\{\itT:\itH\to \itF: \hbox{$T$ linear and continuous}\}$ and denote as
 $\| \itT\|_{\itB}=\displaystyle\sup_{\| u\|\leq 1}\|\itT(u)\|_{\itF} $.
The following assumptions will be required.
\begin{enumerate}
\item[\textbf{A.1}] $\sqrt n (\wmu_n-\mu)=O_\prob(1)$
\item[\textbf{A.2}] $\esp\left[\| X-\mu\|^{-3/2}\right]<\infty$.
\end{enumerate}

\noi \textbf{Theorem \ref{asym}.2.} \textsl{ Under assumptions  \bf A.1 \rm   and  \bf A.2\rm,   we have that
$$\sqrt n(\wGamma^{\Signo}(\wmu_n)-\wGamma^{\Signo}(\mu))=\sqrt n  G_X(\wmu_n-\mu)+o_\prob(1), $$
where $G_X=2F_X-2S_X\in \itB$, with  $F_X$ and $S_X$ defined as follows
\begin{eqnarray}
F_X(u)&=&\esp\left[\frac{\langle(X-\mu),u\rangle}{\| X-\mu\|^4}\; (X-\mu)\otimes (X-\mu) \right]\label{Fx}\\
S_X(u)&=&\frac{1}{2}\left\{u\otimes \esp\left[\frac{X-\mu}{\| X-\mu\|^2}\right]+\esp\left[\frac{X-\mu}{\| X-\mu\|^2}
\right]\otimes u\right\}.\label{Sx}
\end{eqnarray}}

\noi \textbf{Remark \ref{asym}.2.} Assumption  \textbf{A.1} is satisfied for  the spatial median $\mu$, taking $\wmu_n$
as  the averaged of the stochastic gradient algorithm estimator, presented in Cardot \textsl{et al.} (2013), where
the asymptotic distribution of this estimator is obtained (Theorem 3.4).
Regarding  the assumptions  $\esp\left[\| X-\mu\|^{-1}\right]<\infty$ and    \textbf{A.2}, as noted in D\"urre \textsl{et al.} (2014) in the multivariate case, they require that the probability mass is not too strongly concentrated near $\mu$. In particular, assume that the process $X$ has a finite Karhunen--Lo\`eve expansion, $X=\mu+ \sum_{k=1}^q y_k \phi_k$ where $ \phi_k\in \itH$ are orthonormal and $y_k$ are random variables, then $E_{\nu}=\esp\left[\| X-\mu\|^{-\nu}\right]=\esp\left[\| \by\|^{-\nu}\right]$, with $\by=(y_1, \dots, y_q)\trasp$. Hence $E_\nu<\infty$ for $\nu=1, 3/2$  when $\by$ has a bounded density at $0$ while a weaker requirement may be given when $\by$ has an elliptical distribution (see Remark V  in  D\"urre \textsl{et al.}, 2014). For properly infinite--dimensional processes, $E_\nu<\infty$ if  there exists an orthonormal basis $\{\Psi_k\}_{k\ge 1}$ in $\itH$ such that, for some $q\ge 1$, the random vector $\by=(\langle X, \Psi_1\rangle, \dots, \langle X, \Psi_q\rangle)\trasp$ is such that $\esp\left[\| \by\|^{-\nu}\right]<\infty$.   For elliptical   distributed random elements, one may take as basis $\{\Psi_k\}_{k\ge 1}$  the   eigenfunctions of the scatter operator defining the distribution. When the scatter operator of the elliptical distribution has not a finite rank, we have that    $X=\mu\,+\, V\,Y$, where $Y$ is a zero mean Gaussian random element   with covariance operator $\Gamma$ and  the random variable $V>0$ is independent of $Y$. Hence, $E_\nu= \esp V^{-\nu} \, \esp \|Y\|^{-\nu}$ and $E_\nu<\infty$ if and only if $\esp V^{-\nu}<\infty$. In particular, when $V$ is such that $k/V^2\sim \chi_k^2$, which corresponds to the functional version of a multivariate $\itT-$distribution with $k$ degrees of freedom, we have that $E_\nu<\infty$. 

\vskip0.1in

\noi \textbf{Remark \ref{asym}.3.}  It is worth noticing that when  $F_X\equiv 0$ and $\esp\left[({X-\mu})\,{\| X-\mu\|^{-2}}\right]=0$, Theorem \ref{asym}.2  provides an extension  to  the functional data setting of  the result given in Theorem 2 of  D\"urre \textsl{et al.} (2014). More precisely, in this case $\sqrt n(\wGamma^{\Signo}(\wmu_n)-\wGamma^{\Signo}(\mu))=o_\prob(1)$  meaning that the asymptotic behaviour of the spatial sign covariance  operator  is not affected by the imputation of a location estimator.  In particular, if $X$ has a symmetric distribution around its spatial median, meaning that $X-\mu$ and $\mu-X$ have the same distribution and $\wmu_n$ stands for the estimator defined in  Cardot \textsl{et al.} (2013), then $\sqrt n(\wGamma^{\Signo}(\wmu_n)-\wGamma^{\Signo}(\mu))=o_\prob(1)$, so that the asymptotic distribution of $\wGamma^{\Signo}(\wmu_n)$ can be obtained from that of $\wGamma^{\Signo}(\mu)$ using the Central Limit Theorem. 

In particular, for elliptical families this assertion holds.  Furthermore, if $X \, \sim \,\itE(\mu, \Gamma,\varphi)$ and   $\Gamma$ has not a finite rank, using that $X$ can be written as   $X=\mu\,+\, V\,Y$, where $Y$ is   a zero mean Gaussian random element $Y$ with covariance operator $\Gamma$ and   $V>0$ is a random variable independent of $Y$, we get that
  $\Gamma^{\Signo}(\mu)= \esp[s(Y) \otimes s(Y)]$  the sign operator of the process $Y$. Furthermore, noticing that 
$$\wGamma_n^{\Signo}(\mu)=\frac{1}{n}\sum_{i=1}^n s(Y_i) \otimes s(Y_i)\,,$$
we have that $\sqrt n (\wGamma^{\Signo}(\wmu_n)-\Gamma^{\Signo}(\mu))$ converges in $\itF$  to a  zero mean Gaussian element with  covariance
operator equal to the covariance operator of  $s(Y) \otimes s(Y) $ if  the estimator of the location $\mu$ is the functional median $\wmu_n$  defined in  Cardot \textsl{et al.} (2013). 

\vskip0.1in
For multivariate data,  Theorem 2 of  D\"urre \textsl{et al.} (2014) gives the asymptotic distribution of the spatial sign operator. Corollary  \ref{asym}.2 below extends this results to the functional setting. In the general situation in which one cannot guarantee that  $F_X\equiv 0$ and $\esp\left[({X-\mu})\,{\| X-\mu\|^{-2}}\right]=0$, a joint asymptotic distribution between the location parameter estimator and $\wGamma^{\Signo}(\mu)$ is needed.

\noi \textbf{Corollary  \ref{asym}.1.} \textsl{Assume that \textbf{A.2} holds and that
$\left(\sqrt n \left(\wmu_n-\mu\right), \sqrt n\left(\wGamma^{\Signo}(\mu)-\Gamma^{\Signo}(\mu)\right)\right)\convdist Z$, where  $Z$ is a zero mean Gaussian 
random object in $\itH\times\itF$, with covariance operator
$\Upsilon:\itH\times \itF\to \itH\times \itF.$
Then,  $\sqrt n (\wGamma^{\Signo}(\wmu_n)-\Gamma^{\Signo}(\mu))$ converges in $\itF$  to a  zero mean Gaussian element with  covariance
operator given by $(G_X\Pi_1+\Pi_2)\Upsilon (G_X\Pi_1+\Pi_2)^*$,
where $\Pi_i$, for $i=1,2$,  are the projection operators from $\itH\times\itF$ to $\itH$ and $\itF$,  respectively.
Moreover,  $G_X^*=2F_X^*-2S_X^*$ with  $F_X^*$ and $S_X^*$ the adjoint operators of $F_X$ and $S_X$, respectively given by
\begin{eqnarray*}
S_X^*(\Upsilon)&=& \frac 12 \left\{\Upsilon\left(\esp\left[\frac{X-\mu}{\| X-\mu\|^2}\right]\right)+
\Upsilon^*\left(\esp\left[\frac{X-\mu}{\| X-\mu\|^2}\right]\right)\right\}\\
F_X^*(\Upsilon)&=& \esp\left[\frac{ \langle (X-\mu)\otimes (X-\mu),\Upsilon\rangle_{\itF} }{\| (X-\mu)\|^4}\;  (X-\mu)\right].
\end{eqnarray*}
}

\section{Applications}  \label{aplication}
In this section, we consider two applications of the results obtained in Section \ref{asym}. The first one is a result allowing to derive the asymptotic behaviour of the principal direction estimators obtained as the eigenfunctions of $\wGamma^{\Signo}(\wmu_n)$. The second one uses the asymptotic distribution of the sample spatial sign operator to obtain a test for equality among sign covariance operators.

\subsection{On the asymptotic behaviour of the spherical  principal direction estimators}
Robust estimators of the principal directions for functional data have been extensively studied since the spherical principal components proposed in Locantore \textsl{et al.} (1999) and  studied in Gervini (2008). As mentioned in the Introduction, Fisher--consistency of several proposals including the spherical principal directions has been studied in a framework more general than Gaussian random elements, without requiring finite moments, such as that given by elliptically distributed random elements. 

When considering the spherical principal directions two possible situations may arise: either (a)
the distribution is concentrated on a finite--dimensional subspace or (b) the rank of $\Gamma^{\Signo}(\mu)$ is infinite, where $\mu$ stands for the location parameter of $X$  which is typically the functional median. Gervini (2008) showed that the spherical principal direction estimators are Fisher--consistent for the principal directions when the process admits a Karhunen--Lo\`eve
expansion with only finitely many terms, while Boente \textsl{et al.} (2014) derived that the spherical principal components are in fact Fisher--consistent
for any elliptical distribution. More precisely, assume that either: 
\begin{itemize}
\item[a)] $X=\mu +\sum_{k=1}^q \lambda_k^{1/2}\xi_k \phi_k$, where  $\lambda_1\ge   \dots\ge \lambda_k>0$, $\phi_k\in \itH$ are orthonormal and $\xi_k$ are random variables such that $(\xi_1, \dots, \xi_q)\trasp$ has symmetric and exchangeable marginals, 
\item[b)] $X\sim \itE(\mu, \Gamma,\varphi)$ and denote $\lambda_1\ge \lambda_2\ge \dots$  the eigenvalues of the scatter operator $\Gamma$ with associated eigenfunctions $\phi_j$, 
\end{itemize}
hold. 
Note that in the situation a), the scatter operator $\Gamma=\sum_{i=1}^q \lambda_k \phi_k$ has finite rank. As shown in  Gervini (2008)  and  Boente \textsl{et al.} (2014),    the eigenfunctions of $\Gamma^{\Signo}(\mu)$ are those of $\Gamma$ and in the same order. More precisely, if  $\lambda_1^{\Signo}\ge \lambda_2^{\Signo}\ge \dots$ stand for the ordered eigenvalues of $\Gamma^{\Signo}(\mu)$, under a) or b), we have that $\phi_k$ is the eigenfunction of $\Gamma^{\Signo}(\mu)$ related to the eigenvalue  $\lambda_k^{\Signo}$,  meaning that the spatial principal directions are Fisher--consistent. Moreover, we also have that $\lambda_j^{\Signo}> \lambda_{j+1}^{\Signo}$ if $\lambda_j > \lambda_{j+1}$.

Beyond Fisher--consistency, consistency and order of consistency are also desirable properties for any robust procedure. However, for most of the proposed methods only consistency results were obtained. In this section, we derive the asymptotic distribution of the spherical principal direction estimators,  which correspond to the eigenfunctions of the spatial sign operator estimator. In this sense, our result  provides the first asymptotic normality result for robust  principal direction estimators in a general setting.

Even though the asymptotic behaviour of the eigenfunctions of $\wGamma^{\Signo}(\mu)$ can easily be obtained from the central limit theorem and the results in Dauxois \textsl{et al.} (1982), Theorem \ref{asym}.3 states that this asymptotic behaviour may not be the same when location is unknown and estimated. However, it should be noticed that for elliptical distributed random elements or under the symmetry assumptions required in Gervini (2008) to ensure Fisher consistency, we have that the asymptotic behaviour of the eigenfunctions of $\wGamma^{\Signo}(\wmu_n)$ is that of  the eigenfunctions of $\wGamma^{\Signo}(\mu)$, since as mentioned in Remark \ref{asym}.3 $\sqrt n(\wGamma^{\Signo}(\wmu_n)-\wGamma^{\Signo}(\mu))=o_\prob(1)$. 

 Similar arguments to those considered in Dauxois \textsl{et al.} (1982)  and Corollary  \ref{asym}.1  allow to obtain the asymptotic distribution of the   spatial principal direction estimators  not only for elliptical families. For that purpose, denote $\{\lambda_j^{\Signo}\}_{j\ge 1}$  the sequence of eigenvalues of $\Gamma^{\Signo}(\mu)$ ordered in decreasing order and as  $\{\phi_j^{\Signo}\}_{j\ge 1}$ their related eigenfunctions. Let  $\wphi_1^{\Signo},  \wphi_2^{\Signo}, \dots$ be the eigenfunctions of $\wGamma^{\Signo}(\wmu_n)$ related to the ordered eigenvalues $\wlam_1^{\Signo}\ge \wlam_2^{\Signo}\ge \dots$. Recall that if the process has an elliptical distribution with scatter operator $\Gamma$, $\phi_j^{\Signo}=\phi_j$ the $j-$th eigenfunction  of $\Gamma$.

Define $\Lambda_i=\{j\in \natu: \lambda_j^{\Signo} =\lambda_i^{\Signo} \}$, $\Lambda=\{i\in \natu: \mbox{card}(\Lambda_i)=1 \}$  and the projection operators $\Pi_i^{\Signo}=\sum_{j\in \Lambda_i} \phi_j^{\Signo}\otimes \phi_j^{\Signo}$  and $\wPi_i^{\Signo}=\sum_{j\in \Lambda_i} \wphi_j^{\Signo}\otimes \wphi_j^{\Signo}$. The following result is a direct consequence of Propositions 3, 4, 6 and 10 in Dauxois \textsl{et al.} (1982) and  Corollary  \ref{asym}.1.   Taking into account that the $i-$th principal direction is defined up to a sign change when the eigenvalue $\lambda_i^{\Signo}$ has multiplicity one, in the sequel,  we choose the direction of the eigenfunction estimator so that $\langle \wphi_i^{\Signo},  \phi_i^{\Signo}\rangle >0$. 
 
\vskip0.1in
\noi \textbf{Proposition \ref{aplication}.1} \textsl{Assume that \textbf{A.2} holds and that
$\left(\sqrt n \left(\wmu_n-\mu\right), \sqrt n\left(\wGamma^{\Signo}(\mu)-\Gamma^{\Signo}(\mu)\right)\right)\convdist Z$, where  $Z$ is a zero mean Gaussian random object in $\itH\times\itF$, with covariance operator $\Upsilon:\itH\times \itF\to \itH\times \itF$. Denote as $\Pi_i$, for $i=1,2$,  the projection operators from $\itH\times\itF$ to $\itH$ and $\itF$,  respectively and as $U$ a zero mean Gaussian random object in $\itF$ with covariance operator $\Upsilon^{\Signo}=(G_X\Pi_1+\Pi_2)\Upsilon (G_X\Pi_1+\Pi_2)^*$. Then, we have that
\begin{enumerate}[label=\alph*)]
\item   $\wPi_i^{\Signo}\convpp \Pi_i^{\Signo}$ in $\itF$. Moreover, for any $i\in \Lambda$, $\wphi_i^{\Signo}\convpp\phi_i^{\Signo}$ in $\itH$.
\item $\sqrt{n}\left(\wPi_i^{\Signo}-\Pi_i^{\Signo}\right)$ converges in distribution to the zero mean Gaussian random element of $\itF$ given by $ \Delta_i U \Pi_i^{\Signo}+ \Pi_i^{\Signo} U \Delta_i$ where
$$\Delta_i= \sum_{\ell \in \Lambda-\Lambda_i} \frac{1}{\lambda_i^{\Signo}-\lambda_\ell^{\Signo}} \phi_\ell^{\Signo} \otimes \phi_\ell^{\Signo}\,.$$
Furthermore, when $i\in \Lambda$, we have that   $\sqrt{n}\left(\wphi_i^{\Signo}-\phi_i^{\Signo}\right)\convdist (\Delta_i U)(\phi_i)$, which is a zero mean Gaussian process in $\itH$.
\end{enumerate}
}

Note that when $ i\in \Lambda$, $\Delta_i= \sum_{\ell \ne  i} \left\{1/\left(\lambda_i^{\Signo}-\lambda_\ell^{\Signo}\right)\right\} \phi_\ell^{\Signo} \otimes \phi_\ell^{\Signo} $.

\subsection{Tests for equality of the sign covariance operators}{\label{test}}
The asymptotic distribution of the spatial covariance operator given in Corollary  \ref{asym}.1 allows to construct a test for equality between
spatial  covariance  operators between two different
populations. More precisely, assume that we have independent observations
$X_{i,1},\cdots,X_{i,n_i}$, $i=1, 2$ such that $X_{i,j}\sim X_i$, $1\le j\le n_i$ with location parameter $\mu_i$. For the sake of simplicity, let us denote $\Gamma^{\Signo}_i=\esp[s(X_{i}-\mu_i)\otimes s(X_{i}-\mu_i)]$  the  spatial sign covariance  operator of
the $i$-th population.
We are interested in  testing  the null hypothesis
\begin{equation}
H_0:\Gamma^{\Signo}_1=\Gamma^{\Signo}_2\quad\mbox{against}\quad H_1:  \Gamma^{\Signo}_1\neq\Gamma^{\Signo}_2\;\;.
\label{testk}
\end{equation}
As in Boente \textsl{et al.}(2017), we will reject the null hypothesis when the difference  between the estimated spatial sign covariance  operators is large. Namely, if $\wGamma^{\Signo}_i$ stands for a consistent
estimator of $\Gamma^{\Signo}_i$ based on $X_{i,1},\cdots,X_{i,n_i}$,  $i=1,2$, we define
\begin{equation}
T_{n}^{\Signo}=n  \|\wGamma^{\Signo}_2-\wGamma^{\Signo}_1\|_{\itF}^2\,,
\label{tkn}
\end{equation}
where $n=n_1+n_2$. The asymptotic distribution of $T_{n}^{\Signo}$ can be obtained from the asymptotic distribution of
$\sqrt n(\wGamma^{\Signo}_i-\Gamma^{\Signo}_i)$, as stated in the following proposition, which can be considered as a robust version of Corollary 1 in Boente \textsl{et al.} (2017). Its proof can be obtained using Theorem 1 from the above--mentioned  paper.

\vskip0.1in
\noi \textbf{Proposition \ref{aplication}.2} \textsl{
Let $X_{i,1} ,\cdots,X_{i,n_i} \in \itH$,   $i=1, 2$, be independent observations from two independent populations  with location parameter $\mu_i$ and spatial sign covariance  operator $\Gamma^{\Signo}_i$.  Assume that ${n_i}/n\to \tau_i$ with $\tau_i\in (0,1)$
where $n=n_1+n_2$. Let $\wGamma^{\Signo}_i$ be  independent   estimators of  the $i-$th population spatial sign covariance  operator
such that $\sqrt n_i\left(\wGamma^{\Signo}_i-\Gamma^{\Signo}_i\right)\convdist  U_i$,  with $ U_i$ a zero mean Gaussian random element with
covariance operator $\Upsilon_i$. Denote $ \Upsilon_{\weigh}:{\itF}\to {\itF}$ the linear operator defined as $
\Upsilon_{\weigh}= ({1}/{\tau_1})\Upsilon_1 +({1}/{\tau_2})\Upsilon_{2}$ and let $\{\theta_\ell\}_{\ell\,\ge 1} $ stand for the sequence of eigenvalues of $  \Upsilon_{\weigh}$ ordered in decreasing order. 
\new Then, we have that
$n \|(\wGamma^{\Signo}_2-\Gamma^{\Signo}_2)-(\wGamma^{\Signo}_1-\Gamma^{\Signo}_1)\|_\itF^2 \convdist  \sum_{\ell\,\geq 1}\theta_\ell  Z_\ell^2$, with $Z_\ell\sim N(0,1)$  independent. In particular, if $H_0: \Gamma^{\Signo}_1=\Gamma^{\Signo}_2$ holds, we have that
\begin{equation}
n \| \wGamma^{\Signo}_2 -\wGamma^{\Signo}_1 \|_\itF^2\convdist  \sum_{\ell\,\geq 1}\theta_\ell  Z_\ell^2.
\label{disti}
\end{equation}
}

The asymptotic results obtained in Section \ref{asym}, in particular Corollary \ref{asym}.1 and Remark \ref{asym}.3,  invite to consider as estimators of the sign operator $\wGamma^{\Signo}_i=({1}/{n_i})\sum_{j=1}^{n_i} s(X_{i,j}- \wmu_{n_i}) \otimes s(X_{i,j}-  \wmu_{n_i}) $,
with $\wmu_{n_i}$ any consistent estimators of the functional median $\mu_i$ of the process $X_{i}$  satisfying \textbf{A.1}, for instance, the spatial median are given in  Cardot \textsl{et al.} (2013) (see  Remark \ref{asym}.2).
In such a case, as noted in Boente \textsl{et al.} (2017), equation (\ref{disti}) motivates the use of the bootstrap methods, to decide whether to reject the null hypotheses, as follows:
\begin{itemize}
\item[] \textbf{Step  1.} For $i=1,2$ and given the sample $X_{i,1},\cdots,X_{i,n_i}$, let  $\wUps_{i}$ be consistent estimators of $\Upsilon_i$. Define
$\wUps_{\weigh}=\wtau_1^{\,-1}\wUps_1+\wtau_2^{\,-1}\wUps_2$ with $\wtau_i=n_i/({n_1+n_2})$.
\item[]\textbf{Step  2.} For $1\le \ell \le q_n$ denote by  $\wtheta_\ell $  the positive eigenvalues of  $\wUps_{\weigh}$.
\item[]\textbf{Step  3.} Generate $Z^*_1,\dots,Z^*_{q_n}$ i.i.d. such that $Z^*_i \sim N(0,1)$ and let ${\itU}^*_{n}=\sum_{j=1}^{q_n}\wtheta_j {Z^*_{j}}^2$.
\item[]\textbf{Step  4.}  Repeat \textbf{Step  3}  $N_{\boot}$ times, to get $N_{\boot}$ values of ${\itU}_{nr}^*$ for $1\leq r\leq N_{\boot}$.
\end{itemize}
The $(1-\alpha)-$quantile of the asymptotic null distribution of $T_{n}^{\Signo}$ can be approximated by the $(1-\alpha)-$quantile of
the empirical distribution of ${\itU}_{nr}^*$ for $1\leq r\leq N_{\boot}$.  Besides, the $p-$value can be estimated by $\widehat{p}=s/{N_{\boot}}$ where $s$ equals the number of  ${\itU}^*_{nr}$ which are larger or equal than the observed value of the statistic $T_{n}^{\Signo}$. The validity of the bootstrap procedure can be derived from Theorem 3 in  Boente \textsl{et al.} (2017) if the estimators of estimators of $\Upsilon_{\weigh}$ are such that $\sqrt{n}\|\wUps_{\weigh}-\Upsilon_{\weigh}\|_\itF=O_\prob(1)$ ensuring that the asymptotic significance level of the test based on the bootstrap critical value is indeed $\alpha$. A possible choice for  $\wUps_{i}$, $i=1,2$ is the sample covariance operator of $Y_i=s(X_{i}-\mu_i)\otimes s(X_{i}-\mu_i)$.

\vskip0.2in

\noi \textbf{Remark \ref{aplication}.1.} Proposition \ref{aplication}.2 ensures that, under mild assumptions, it is possible to provide a test to decide if $\Gamma^{\Signo}_1=\Gamma^{\Signo}_2$. An important point to highlight is what this null hypothesis represents, for instance, in terms of the covariance   operators of the two populations, when they exist.  Let us consider the situation in which the two populations have an elliptical distribution, that is, $X_i\sim \itE(\mu_i, \Gamma_i, \varphi_i)$, for $i=1,2$. Recall that the eigenfunctions of $\Gamma^{\Signo}_i$ are those of $\Gamma_i$ and in the same order, while  the eigenvalues of the sign covariance  operator $\Gamma^{\Signo}_i$, denoted $\lambda_{i,\ell}^{\Signo}$, are  shrunk with respect to those of $\Gamma_i$ (that are denoted as $\lambda_{i,\ell}$)  as follows
 \begin{equation}
  \lambda_{i,\ell}^{\Signo}  = \lambda_{i,\ell} \, \esp\left(\frac{\xi^2_{i,\ell}}{\sum_{j \geq 1} \lambda_{i,j} \xi_{i,j}^2}\right)\,,
  \label{lambdaiS}
  \end{equation}
where $\xi_{i,j}=\lambda_{i,\ell}^{-1/2}\,\langle X_{i}-\mu_i, \phi_{i,j}\rangle$ with $\phi_{i,\ell}$ the eigenfunction of $\Gamma_i$. 

Assume that $\varphi_1=\varphi_2$, that is, if the two populations have the same underlying distribution up to location and scatter. Note that,  if the scatter operators are proportional, i.e., if $\Gamma_2=\rho\, \Gamma_1$  for some positive constant $\rho$, then $\Gamma^{\Signo}_1=\Gamma^{\Signo}_2$.  Thus, when the two populations have the same elliptical distribution up to changes in location and scatter, the test based on $T_{n}^{\Signo}$  provides a way for testing proportionality of the scatter operators, even when second moments do not exist.   It is worth noticing that, when second moment exists the covariance operator of $X_i$ is up to a constant equal to $\Gamma_i$, hence the statistic  $T_{n}^{\Signo}$  allows to test proportionality between the two covariance operators. Note that when $\Gamma^{\Signo}_1=\Gamma^{\Signo}_2$, both scatter operators have the same rank and share the same eigenfunctions. Furthermore, if the scatter operators have finite rank, from Proposition 1 in D\"urre \textsl{et al.} (2016), we get that $\Gamma^{\Signo}_1=\Gamma^{\Signo}_2$ if and only if $\Gamma_2=\rho\, \Gamma_1$  for some positive constant $\rho$. Hence, for finite rank scatter operators,  testing proportionality of the scatter operators is equivalent to testing equality of the spatial sign operators.

\subsubsection{Monte Carlo study}\label{simu}
This section contains the results of a simulation study devoted to illustrate the finite--sample performance of the test procedure described in Section  \ref{test}, under the null hypothesis and  different alternatives, when atypical data are introduced in the samples.   The numerical study also aim to compare the performance of the sign operator testing procedure   with that based on the sample covariance operator   introduced in Boente \textsl{et al.} (2017).

We have performed $N = 1000$ replications taking samples of size $n_i = 100$, $i=1, 2$. The generated samples $X_{i,1},\cdots,X_{i,n_i}$, $i=1, 2$ are such that $X_{i,j}\sim X_i\in L^2(0,1)$.  In all cases,  each trajectory was observed at $m = 100$ equidistant points in the interval $[0, 1]$ and we performed $N_{\boot} =5000$ bootstrap replications. To summarize the tests performance,  we compute the observed frequency of rejections over replications with nominal level $\alpha=0.05$.

\subsubsection*{Simulation settings}
The distribution of the two populations correspond, under the null hypothesis,  to independent centred  Brownian motion processes,  denoted from now on as $\itB\itM(0,1)$.  Hence, both processes have the same spatial sign operators and also the same covariance operators. On the other hand, to check the test power performance, we consider   the same  alternatives as in Boente \textsl{et al.} (2017)  and also Gaussian alternatives. More precisely,  we generate independent observations  $X_{1,j}\sim X_1$, $1\le j\le n_1$, and $X_{2,j}\sim X_2$, $1\le j\le n_2$,  such that  $X_1$ has the distribution of a centred Brownian motion denoted $\itB\itM (0,1)$ while the second population has a distribution according to the one of the following models
\begin{itemize}
\item \textbf{Model 1:} $X_2\sim  Y_1+\delta_n \, Y_2^2$, where $Y_1$ and $Y_2$ are independent $Y_i\sim {\mathcal{BM}}(0,1)$, $i=1,2$ and $\delta_n= \Delta\, n^{-1/4}$ with $n=n_1+n_2$ and $\Delta$ takes values from 0 to 8 with step 1 and from 10 to 20 with step 2. The situation $\Delta=0$ corresponds to the null hypothesis in which both processes have a Gaussian distribution. 
\item \textbf{Model 2:} $X_2\sim  Y_1+\delta_n \, Y_2$, where $Y_1$ and $Y_2$ are independent $Y_1\sim {\mathcal{BM}}(0,1)$, $Y_2$ is a Gaussian process with covariance kernel $\cov(Y_2(t),Y_2(s))=\exp(-\,|s-t|/0.2)$ and $\delta_n= \Delta\, n^{-1/4}$ with $n=n_1+n_2$ and $\Delta\in \{0, 0.5, 1, 1.5, 2,2.5,3,4,5\}$. In this case both processes  have a Gaussian distribution under the null and under the alternative which implies that, for each population,  the spatial sign operator   has the same eigenfunctions as the  covariance operator. Moreover,  the eigenvalues of the   spatial operator and of the covariance operator of the $i-$th population are related   through (\ref{lambdaiS}) with $\xi_{i,\ell}\sim N(0,1)$ independent of each other. 
\end{itemize}

To analyse the behaviour when atypical data are introduced in the sample, for each generated sample, we also consider the following contamination. We first generate two independent samples $V_{i,j}\in \real$, $1\le i\le n_j$ and $j=1,2$, such that $V_{i,j} \sim |\itT_1|$,  where $|\itT_1|$ corresponds to the absolute value of an univariate  $\itT$-Student distribution with 1 degree of freedom. We then generated the contaminated samples, denoted $X_{i,j}^{(c)}$,   as follows  $X_{i,j}^{(c)}=(1-B_i)\;X_{i,j}+B_i\;V_{i,j}\;X_{i,j}$ where  $B_i \sim Bi(1,0.1)$ are independent and independent of $(X_{i,j}, V_{i,j})$. Note that under the null hypothesis, both populations have the same elliptical distribution since they can be written as $ W_{i,j}\;X_{i,j} $ with $W_{i,j}=(1-B_i)+B_i\;V_{i,j}$ a positive random variable independent of $X_{i,j}$ and $W_{1,j}\sim W_{2,j}$.

\subsubsection*{The test statistics}
 
We computed two test statistics, the statistic based on the spatial sign operator defined above and the procedure defined in Boente \textsl{et al.} (2017). The test statistic given in this last paper is defined as $T_n=n\,\|\wGamma_1-\wGamma_2\|^2$ where $\wGamma_i=(1/n_i) \sum_{j=1}^{n_i} (X_{i,j}-\overline{X}_{i})\otimes (X_{i,j}-\overline{X}_{i})$ is the sample covariance operator. This testing method is designed to test equality of the two populations covariance operators, which is fulfilled when $\Delta=0$. On the other hand, the statistic $T_n^{\Signo}$  defined in (\ref{tkn}) is designed to test  equality of the spatial operators, that is, $\Gamma_1^{\Signo}=\Gamma_2^{\Signo}$. As mentioned in Remark \ref{aplication}.1, this null hypothesis is fulfilled when the scatter operators related to the elliptical distribution are proportional which holds when $\Delta=0$, both for clean and contaminated samples. When computing the spatial sign operators $\wbGa_i^{\Signo}$, we center the data with the functional median computed through the function \textsl{l1median} from the R package \textsl{pcaPP}.

The testing procedure requires bootstrap calibration. For that purpose, following the procedure described in Boente \textsl{et al.} (2017),  we  project the centred  data onto the   $M$  largest principal components of the pooled  operators $\wGamma_{\pool}$, where the pooled operator was adapted to the testing procedure used. More precisely,   $\wGamma_{\pool}=n^{-1}\sum n_i \wbGa_i$  when the test statistic is based on the sample covariance matrices, while  $\wGamma_{\pool}=n^{-1}\sum n_i \wbGa_i^{\Signo}$, when the test statistic corresponds to  the sample sign operator. The covariance operator of each  estimator, denoted $\Upsilon_{\weigh}$ in Proposition  \ref{aplication}.2 for the spatial operator, is then estimated through a finite dimensional matrix. We choose different values of the number of principal directions $M=  3, 10, 20$ and $30$ to study the dependence on the finite--dimensional approximation considered. As noted in Boente \textsl{et al.} (2017), the value $q_n$ used in \textbf{Step 2} equals $q_n=M(M+1)/2$. With the selected number of principal directions, we explained more than 80\% of the total variability (see Table \ref{tab:percent}).

\begin{table}[ht!]
    \centering
   \small
\begin{tabular}{c|c|cccc|cccc|}
 \hline
 & & \multicolumn{4}{c|}{Clean samples} & \multicolumn{4}{c|}{Contaminated samples}\\\hline
 & $\Delta$& \multicolumn{4}{c|}{$M$}   & \multicolumn{4}{c|}{$M$}  \\\hline
 & & 3 & 10 & 20 & 30 & 3 & 10 & 20 & 30   \\
  \hline
$T_{\boot,M}$ & 0 & 0.934 & 0.981 & 0.991 & 0.995 & 0.962 & 0.992 & 0.996 & 0.998\\
$T^{\Signo}_{\boot,M}$ & 0 &  0.828 & 0.950 & 0.979 & 0.989 & 0.828 & 0.950 & 0.979 & 0.989 \\ 
   \hline
\end{tabular}
\caption{\small \label{tab:percent} Percentage of the total variance explained by the first $M$ principal components when using the test $T_{\boot,M}$ or $T^{\Signo}_{\boot,M}$.}
\end{table}

When  the populations have a Gaussian distribution, the asymptotic covariance operator of the sample covariance operator $\wGamma_i=(1/n_i) \sum_{j=1}^{n_i} (X_{i,j}-\overline{X}_{i})\otimes (X_{i,j}-\overline{X}_{i})$ can be  estimated using the eigenvalues and eigenfunctions of  $\wGamma_i$.  Taking into account that, under the null hypothesis, the processes are Gaussian, we have also used this approximation when considering the sample covariance operator.

From now on we denote as $T_{\boot,M}$ and $T^{\Signo}_{\boot,M}$, for $M=3,10,20$ and $30$ the bootstrap calibration of the statistics $T_{n}$ and $T_n^{\Signo}$, respectively, computed using $M$ principal components. Finally,  $T_{\boot,\gaus}$ stands for the bootstrap calibration of $T_{n}$ computed using the Gaussian approximation.

\subsubsection*{Simulation results}
For the alternatives given through Models 1 and 2, Tables   \ref{tab:clasica} and \ref{tab:clasicamod2} summarize, respectively,   the frequency of rejection  for the    procedure based on the sample covariance operator for the uncontaminated samples and for the contaminated samples, while those corresponding to the   test based on the sample spatial sign operator are reported in Table  \ref{tab:robusta} and \ref{tab:robustamod2}.

\begin{table}[ht!]
 \small
  \centering
   \setlength{\tabcolsep}{3pt} 
    \renewcommand{\arraystretch}{1.1} 

\begin{tabular}{c|ccccccccccccccc|}
  \hline
$\Delta$ & 0 & 1 & 2 & 3 & 4 & 5 & 6 & 7 & 8 & 10 & 12 & 14 & 16 & 18 & 20 \\ 
  \hline
  & \multicolumn{15}{c|}{Clean samples}\\\hline
$T_{\boot,3}$ & 0.066 & 0.083 & 0.315 & 0.694 & 0.895 & 0.948 & 0.959 & 0.967 & 0.972 & 0.973 & 0.974 & 0.974 & 0.974 & 0.975 & 0.975\\ 
 $T_{\boot,10}$ & 0.065 & 0.082 & 0.299 & 0.681 & 0.890 & 0.942 & 0.957 & 0.962 & 0.969 & 0.971 & 0.972 & 0.973 & 0.973 & 0.973 & 0.973 \\ 
  $T_{\boot,20}$ &  0.061 & 0.081 & 0.296 & 0.671 & 0.885 & 0.941 & 0.956 & 0.961 & 0.965 & 0.968 & 0.971 & 0.973 & 0.973 & 0.973 & 0.971 \\ 
 $T_{\boot,30}$ & 0.060 & 0.079 & 0.290 & 0.666 & 0.882 & 0.940 & 0.956 & 0.961 & 0.964 & 0.967 & 0.970 & 0.973 & 0.973 & 0.971 & 0.971\\ 
  \hdashline
 $T_{\boot,\gaus}$ & 0.050 & 0.064 & 0.333 & 0.801 & 0.975 & 0.998 & 1.000 & 1.000 & 1.000 & 1.000 & 1.000 & 1.000 & 1.000 & 1.000 & 1.000  \\
 \hline
  & \multicolumn{15}{c|}{Contaminated samples} \\\hline
  $T_{\boot,3}$ &  0.011 & 0.016 & 0.037 & 0.082 & 0.152 & 0.217 & 0.268 & 0.291 & 0.315 & 0.364 & 0.394 & 0.412 & 0.436 & 0.451 & 0.461 \\ 
  $T_{\boot,10}$& 0.010 & 0.013 & 0.033 & 0.076 & 0.148 & 0.214 & 0.260 & 0.283 & 0.306 & 0.355 & 0.389 & 0.402 & 0.423 & 0.441 & 0.451 \\ 
 $T_{\boot,20}$ &  0.009 & 0.011 & 0.032 & 0.074 & 0.145 & 0.212 & 0.256 & 0.280 & 0.302 & 0.348 & 0.381 & 0.396 & 0.420 & 0.436 & 0.446\\ 
 $T_{\boot,30}$& 0.009 & 0.011 & 0.031 & 0.074 & 0.144 & 0.209 & 0.253 & 0.279 & 0.300 & 0.345 & 0.379 & 0.395 & 0.417 & 0.434 & 0.445\\   \hdashline
 $T_{\boot,\gaus}$ & 0.843 & 0.836 & 0.856 & 0.923 & 0.961 & 0.974 & 0.979 & 0.984 & 0.989 & 0.994 & 0.997 & 0.996 & 0.998 & 0.999 & 0.999\\
   \hline
\end{tabular} 
\caption{\small \label{tab:clasica}   Frequency of rejection for the bootstrap test  $T_{\boot,M}$ based on the sample covariance operators,  under Model 1, when   $M=3,10,20$ and $30$ principal components are used. The  row labelled $T_{\boot,\gaus}$ reports the frequencies obtained when the eigenvalues $\theta_\ell$ are estimated using that the processes are Gaussian.}
 \end{table}

\begin{table}[ht!]
 \small
  \centering
   \setlength{\tabcolsep}{3pt} 
    \renewcommand{\arraystretch}{1.1} 

\begin{tabular}{c|ccccccccccccccc|}
  \hline
$\Delta$ & 0 & 1 & 2 & 3 & 4 & 5 & 6 & 7 & 8 & 10 & 12 & 14 & 16 & 18 & 20 \\ 
  \hline
  & \multicolumn{15}{c|}{Clean samples}\\\hline
$T_{\boot,3}^{\Signo}$ & 0.046 & 0.046 & 0.084 & 0.165 & 0.262 & 0.375 & 0.487 & 0.588 & 0.674 & 0.776 & 0.833 & 0.872 & 0.901 & 0.925 & 0.932 \\ 
 $T_{\boot,10}^{\Signo}$ & 0.047 & 0.058 & 0.092 & 0.179 & 0.290 & 0.396 & 0.507 & 0.618 & 0.700 & 0.798 & 0.854 & 0.892 & 0.921 & 0.934 & 0.942 \\ 
  $T_{\boot,20}^{\Signo}$ & 0.047 & 0.059 & 0.094 & 0.182 & 0.294 & 0.401 & 0.514 & 0.621 & 0.702 & 0.804 & 0.855 & 0.895 & 0.922 & 0.937 & 0.942 \\ 
 $T_{\boot,30}^{\Signo}$ & 0.048 & 0.058 & 0.094 & 0.183 & 0.294 & 0.402 & 0.515 & 0.622 & 0.703 & 0.806 & 0.856 & 0.895 & 0.922 & 0.938 & 0.944 \\ 
 \hline
  & \multicolumn{15}{c|}{Contaminated samples} \\\hline
  $T_{\boot,3}^{\Signo}$ & 0.048 & 0.047 & 0.076 & 0.127 & 0.218 & 0.327 & 0.448 & 0.529 & 0.599 & 0.711 & 0.786 & 0.834 & 0.873 & 0.889 & 0.899\\ 
  $T_{\boot,10}^{\Signo}$& 0.050 & 0.050 & 0.082 & 0.148 & 0.237 & 0.349 & 0.464 & 0.549 & 0.612 & 0.738 & 0.804 & 0.851 & 0.883 & 0.896 & 0.911 \\ 
 $T_{\boot,20}^{\Signo}$ &  0.052 & 0.052 & 0.082 & 0.149 & 0.243 & 0.351 & 0.470 & 0.553 & 0.616 & 0.743 & 0.806 & 0.853 & 0.885 & 0.900 & 0.911\\ 
 $T_{\boot,30}^{\Signo}$& 0.052 & 0.052 & 0.082 & 0.150 & 0.245 & 0.353 & 0.470 & 0.553 & 0.618 & 0.743 & 0.806 & 0.856 & 0.886 & 0.900 & 0.912 \\ 
   \hline
 \end{tabular}
\caption{\small \label{tab:robusta}   Frequency of rejection for the bootstrap test  $T_{\boot,M}^{\Signo}$ based on the spatial sign operator,  under Model 1,  when   $M=3,10,20$ and $30$ principal directions are used.}
 \end{table}

\begin{table}[ht!]
 \centering
  \footnotesize
   \setlength{\tabcolsep}{2pt} 
    \renewcommand{\arraystretch}{1.1} 
\begin{tabular}{c|ccccccccc|ccccccccc|}
  \hline
$\Delta$ & 0 & 0.5 & 1 & 1.5 & 2 & 2.5 & 3 & 4 & 5 & 0 & 0.5 & 1 & 1.5 & 2 & 2.5 & 3 & 4 & 5  \\ 
  \hline
  & \multicolumn{9}{c|}{Clean samples} & \multicolumn{9}{c|}{Contaminated samples}\\\hline
$T_{\boot,3}$ & 0.066 &  0.057 & 0.063 & 0.131 & 0.289 & 0.678 & 0.949 & 1.000 & 1.000 & 0.011 &  0.004 & 0.007 & 0.009 & 0.028 & 0.060  & 0.095 & 0.140 & 0.177 \\ 
 $T_{\boot,10}$ & 0.065 & 0.056 & 0.062 & 0.121 & 0.265 &  0.619 & 0.917 & 1.000 & 1.000 & 0.010 &   0.002 & 0.005 & 0.009 & 0.024 &  0.044 & 0.082 & 0.129 & 0.156  \\ 
  $T_{\boot,20}$ &  0.061 & 0.055 & 0.059 & 0.120 & 0.262 & 0.598 & 0.910 & 1.000 & 1.000 & 0.009 & 0.002 & 0.005 & 0.009 & 0.021 &  0.043 & 0.078 & 0.128 & 0.154 \\ 
 $T_{\boot,30}$ & 0.060 & 0.053 & 0.057 & 0.116 & 0.257 & 0.588 & 0.909 & 1.000 & 1.000 & 0.009 &0.002 & 0.005 & 0.009 & 0.021 &  0.041 & 0.076 & 0.128 & 0.153  \\ 
  \hdashline
 $T_{\boot,\gaus}$ & 0.050 &0.039 & 0.048 & 0.099 & 0.228 & 0.559 & 0.904 & 1.000 & 1.000 & 0.843 &0.807 & 0.828 & 0.853 & 0.886 &  0.922 & 0.967 & 0.997 & 0.999\\
 \hline
 \end{tabular} 
\caption{\small \label{tab:clasicamod2}   Frequency of rejection for the bootstrap test  $T_{\boot,M}$ based on the sample covariance operators,  under Model 2, when   $M=3,10,20$ and $30$ principal components are used. The  row labelled $T_{\boot,\gaus}$ reports the frequencies obtained when the eigenvalues $\theta_\ell$ are estimated using that the processes are Gaussian.}
 \end{table}

\begin{table}[ht!]
 \footnotesize
  \centering
    \setlength{\tabcolsep}{2pt} 
    \renewcommand{\arraystretch}{1.1} 
\begin{tabular}{c|ccccccccc|ccccccccc|}
  \hline
$\Delta$ & 0 & 0.5 & 1 & 1.5 & 2 & 2.5 & 3 & 4 & 5 & 0 & 0.5 & 1 & 1.5 & 2 & 2.5 & 3 & 4 & 5  \\ 
  \hline
  & \multicolumn{9}{c|}{Clean samples} & \multicolumn{9}{c|}{Contaminated samples}\\\hline
$T_{\boot,3}^{\Signo}$ & 0.046 & 0.057 & 0.166 & 0.491 & 0.849 &  0.975 & 0.999 & 1.000 & 1.000  & 0.048 & 0.053 & 0.174 & 0.506 & 0.847 & 0.973 & 0.999 & 1.000 & 1.000\\ 
 $T_{\boot,10}^{\Signo}$ & 0.047 & 0.061 & 0.178 & 0.521 & 0.867 & 0.982 & 1.000 & 1.000 & 1.000 & 0.050 &  0.057 & 0.183 & 0.530 & 0.874 &  0.979 & 1.000 & 1.000 & 1.000 \\ 
  $T_{\boot,20}^{\Signo}$ & 0.047 & 0.062 & 0.180 & 0.525 & 0.875 & 0.982 & 1.000 & 1.000 & 1.000 &  0.052 & 0.059 & 0.185 & 0.535 & 0.880 &  0.979 & 1.000 & 1.000 & 1.000 \\ 
 $T_{\boot,30}^{\Signo}$ & 0.048 &  0.063 & 0.180 & 0.526 & 0.875 & 0.982 & 1.000 & 1.000 & 1.000 & 0.052 & 0.059 & 0.187 & 0.535 & 0.881 &  0.979 & 1.000 & 1.000 & 1.000 \\ 
 \hline 
  \end{tabular}
\caption{\small \label{tab:robustamod2}   Frequency of rejection for the bootstrap test  $T_{\boot,M}^{\Signo}$ based on the spatial sign operator,  under Model 2,  when   $M=3,10,20$ and $30$ principal directions are used.}
 \end{table}

As noted in Boente \textsl{et al.} (2017) when using the    Gaussian  approximation the test based on the sample covariance operators shows an improvement in size for uncontaminated samples. However, when contaminating the data the level breaks--down and the test becomes uninformative.

On the other hand, when   projecting the data on the first $M$ principal components,   the empirical size of the   tests based on the bootstrap  calibration either using the sample covariance  or the spatial sign operators is quite close to the nominal one, for uncontaminated samples. To analyse the significance of the empirical size,  we study if the empirical size  is significantly different from the nominal level $\alpha=0.05$ by testing  $H_{0,\pi}: \pi=\alpha$ with nominal level $\gamma$, where $\pi$ stands for the value  such that $\pi_n\convprob \pi$ with $\pi_n$ the empirical size of the considered test. This null hypothesis is rejected at level $\gamma$ versus  $H_{1,\pi}:\pi\ne\alpha$  
if   $\pi_n\notin [a_1(\alpha),a_2(\alpha)]$ where $a_j(\alpha)=\alpha+ (-1)^j z_{\gamma/2}\,\{\alpha(1-\alpha)/N\}^{1/2}$, $j=1,2$. If $H_{0,\pi}: \pi=\alpha=0.05$  is not rejected, the testing procedure based  is considered accurate, while if $\pi_n<a_1(\alpha)$  the testing procedure is  conservative and when $\pi_{n}>a_2(\alpha)$ the test is liberal. For the clean samples, both procedures are  accurate with significance level $\gamma=0.01$. On the other hand, when contaminating the samples, the test based on the sample covariance operator becomes conservative with empirical size not exceeding $0.011$ for any value of $M$, while that based on the sign operator preserves its empirical size.

 Regarding the behaviour under the   alternative,  the procedure based on the spatial sign operator shows a loss of power   with respect to the sample covariance operator when the alternatives follow Model 1. On the other hand, for the Gaussian alternatives, the sign test has a much better performance attaining  a higher   power in particular  when $\Delta$ varies between 1 and 3. For both models, the test $T_{\boot,M}^{\Signo}$ is stable for the considered contaminations, while the procedure based on the sample covariance operator shows an important loss of power, since the frequency of rejection never exceeds $0.5$ or $0.2$ under Models 1 and 2, respectively, for any value of the selected number of principal directions $M$.

\bigskip

\noi\textbf{\small Acknowledgements.} {\small  This research was partially supported by Grants  \textsc{pip}  112-201101-00742 from \textsc{conicet}, \textsc{pict} 2014-0351 and 201-0377 from \textsc{anpcyt} and  20020130100279\textsc{ba} and 20020150200110\textsc{ba}  from the Universidad de
Buenos Aires  at Buenos Aires, Argentina.}

\section*{Appendix}{\label{ape}}
{\setcounter{equation}{0}
\renewcommand{\theequation}{A.\arabic{equation}}
Throughout this section, we will assume that   $\mu=0$, without loss of generality. Furthermore, we will denote as $\wGamma^{\Signo}_0=\wGamma^{\Signo}(0)$, $\Gamma_i( \wmu_n)=s(X_i- \wmu_n)\otimes s(X_i- \wmu_n)$ and $\Gamma_i=s(X_i)\otimes s(X_i)$.

\noi \textsc{Proof of Theorem \ref{asym}.1.}
Note that the strong law of large numbers entails that  it is enough to prove that $\wGamma^{\Signo}( \wmu_n)-\wGamma_0^{\Signo} \convpp 0$. Consider the following random set
$$\itA_n=\{x\in \itH: \| x- \wmu_n\|\geq \frac{1}{2}\| x\|\}.$$
Therefore, we have that
\begin{equation}
\label{endos}
\| \wGamma^{\Signo}( \wmu_n)-\wGamma^{\Signo}_0\|_{\itF}\leq \frac{1}{n}\sum_{X_i\in \itA_n}\|\Gamma_i( \wmu_n)-\Gamma_i\|_{\itF}+
\frac{1}{n}\sum_{X_i\notin \itA_n}\|\Gamma_i( \wmu_n)-\Gamma_i\|_{\itF}=A_{n,1}+A_{n,2}
\end{equation}
To show that $A_{n,1}\convpp 0$, note that straightforward calculations lead to the bound
$$\|\Gamma_i( \wmu_n)-\Gamma_i\|_{\itF}^2=\frac{2}{\| X_i\|^2 \| X_i - \wmu_n\|^2}\left(\| X_i\|^2 \|  \wmu_n\|^2-\langle  \wmu_n,X_i\rangle^2\right)\leq
\frac{4 \| X_i\|^2  \|  \wmu_n\|^2  }{\| X_i\|^2 \| X_i - \wmu_n\|^2}\;.
$$
On the other hand, if  $X_i\in \itA_n$, we have that $\|\Gamma_i( \wmu_n)-\Gamma_i\|_{\itF}^2\leq  {16\|  \wmu_n\| ^2}/{\| X_i\| ^2}$ which implies that
$$\frac{1}{n}\sum_{X_i\in \itA_n}\|\Gamma_i( \wmu_n)-\Gamma_i\|_{\itF}\leq
4\; \| \wmu_n\|  \; \frac{1}{n}\;\sum_{i=1}^n \frac{1}{\| X_i\| }\,.$$
Therefore, using that $\wmu_n\convpp 0$, $\esp\left[\|X\|^{-1}\right]<\infty$ and the strong law of large numbers we conclude that $A_{n,1}\convpp 0$.

In remains to show that  the second term $A_{n,2}$ in the right hand side of (\ref{endos}) converges almost surely to zero. The fact that  $\| \Gamma_i( \wmu_n)\|_\itF=\| \Gamma_i\|_\itF=1$ implies that
$$
A_{n,2}=\frac{1}{n}\sum_{X_i\not\in \itA_n}\|\Gamma_i( \wmu_n)-\Gamma_i\|_{\itF}\leq
 \frac{2}{n} \sum_{i=1}^n Z_{n,i}\,,$$
where $Z_{n,i}=\indica_{\itA^c_n}(X_i)$.

Note that the assumption  $\esp\left[\|X\|^{-1}\right]<\infty$ implies that $\prob(\|X\|=0)=0$.  Hence, for any $\epsilon>0$, let $\delta>0$ be such that $\prob(\| X\|\leq \delta)\leq \varepsilon$ and denote ${Z}_{\delta,i}=\indica_{B_{\delta}}(X_i)$, where  $B_{\delta}=\{\| x\|\leq \delta\}$.
Then,
$$
\frac{1}{n} \sum_{i=1}^n Z_{n,i}\leq\frac{1}{n} \sum_{i=1}^n Z_{\delta,i}+\frac{1}{n} \sum_{i=1}^n (Z_{n,i}-Z_{\delta,i})_+=B_{n,1}+B_{n,2}\;,
$$
where    $a_+=\max(a,0)$. The strong law of large numbers entails that $B_{n,1} \convpp \prob(\| X\|\leq \delta)\leq \varepsilon$. To show that $B_{n,2}\convpp 0$, note that $\{\|  \wmu_n\|\leq \delta/2\}\subset\{(Z_{n,i}-Z_{\delta,i})_+=0\}$. Hence, using that $\wmu_n\convpp 0$, we get that there exists a null probability set $\itN$ such that for $\omega \notin \itN$, there exists $n_0$ such that, for all $n>n_0$, $\|  \wmu_n\|\leq \delta/2$ implying that $B_{n,2}=(1/n) \sum_{i=1}^n (Z_{n,i}-Z_{\delta,i})_+=0$ and concluding the proof. \square

\noi \textsc{Proof of Theorem \ref{asym}.2.}
Note that  $\Gamma_i( \wmu_n)-\Gamma_i$ can be written as follows
\begin{eqnarray}
 \Gamma_i( \wmu_n)-\Gamma_i &=& \| X_i\|^{-2} \left\{\| X_i\|^2\Gamma_i( \wmu_n)-X_i\otimes X_i\right\}\label{gamas}\\
 &=& \| X_i\|^{-2} \left\{\left[\| X_i- \wmu_n\|^2+\|  \wmu_n\|^2+2 \langle X_i- \wmu_n, \wmu_n \rangle \right]\Gamma_i( \wmu_n)-X_i\otimes X_i\right\}\nonumber\\
 &=& \| X_i\|^{-2} \left\{ \wmu_n\otimes  \wmu_n- \wmu_n \otimes X_i-X_i\otimes  \wmu_n+\,\left(2\,\langle X_i, \wmu_n\rangle -\|  \wmu_n\|^2\right)\,\Gamma_i( \wmu_n)\right\}\nonumber\;.
\end{eqnarray}
Therefore, $
\sqrt{n} \left( \wGamma^{\Signo}( \wmu_n) - \wGamma^{\Signo}_0\right)=(1/\sqrt{n})\sum_{i=1}^n\left(\Gamma_i( \wmu_n)-\Gamma_i\right)=S_{n,1}-S_{n,2}-S_{n,3}+2\,S_{n,4}-S_{n,5}$, 
where
\begin{eqnarray*}
S_{n,1} &=& 
\frac{1}{\sqrt n}\sum_{i=1}^n\frac{ \wmu_n\otimes  \wmu_n}{\| X_i\|^{2}}=n\;( \wmu_n\otimes  \wmu_n )\left( \frac{1}{n^{3/2}}\displaystyle \sum_{i=1}^n\frac{1}{\| X_i\|^{2}}\right)\\
S_{n,2} &=& \frac{1}{\sqrt n}\sum_{i=1}^n  \frac{  \wmu_n \otimes X_i}{\| X_i\|^{2}}=\sqrt n \; \wmu_n \otimes \left(\frac{1}{ n}\displaystyle\sum_{i=1}^n \frac{X_i}{\| X_i\|^{2}}\right)\\
S_{n,3}&=&  \frac{1}{\sqrt n}\sum_{i=1}^n \frac{ X_i\otimes  \wmu_n}{\| X_i\|^{2}}=  \frac{1}{ n}\displaystyle\sum_{i=1}^n \frac{X_i}{\| X_i\|^{2}}  \otimes  \sqrt n \; \wmu_n\\
S_{n,4} &=& \frac{1}{\sqrt n} \sum_{i=1}^n \frac{\langle X_i, \wmu_n\rangle }{\| X_i\|^{2}}\Gamma_i( \wmu_n)\\
 S_{n,5} &=&\|  \wmu_n\|^2\frac{1}{\sqrt n}\sum_{i=1}^n\frac{\Gamma_i( \wmu_n)}{\| X_i\|^{2}} = n \|  \wmu_n\|^2 \frac{1}{n^{3/2}}\sum_{i=1}^n\frac{\Gamma_i( \wmu_n)}{\| X_i\|^{2}} \;.
\end{eqnarray*}
Note that \textbf{A.2} entails that $\esp V_i^{2/3} <\infty$ where $V_i=1/\|X_i\|^2$, so the Marcinkiewicz's strong law of large numbers implies that $n^{-3/2} \sum_{i=1}^n 1/{\| X_i\|^{2}}\convpp 0$. Hence, Assumptions \textbf{A.1} and \textbf{A.2} together with the strong law of large numbers and the fact that $\|\Gamma_i( \wmu_n)\|_\itF=1$ entail that  $S_{n,j}\convprob 0$ for $j=1,5$. 

The   decomposition of $\Gamma_i( \wmu_n)-\Gamma_i$ obtained in (\ref{gamas}) entails that $S_{n,4}$ can be written as $S_{n,4}=S_{n41}+S_{n42}-S_{n43}-S_{n44}+S_{n45}-S_{n46}$, where
\begin{eqnarray*}
S_{n41}=\frac{1}{\sqrt n}\sum_{i=1}^n\frac{\langle X_i, \wmu_n\rangle}{\| X_i\|^{4}}\;  X_i\otimes X_i &\qquad & 
S_{n42}=\frac{1}{\sqrt n}\sum_{i=1}^n \frac{\langle X_i, \wmu_n\rangle }{\| X_i\|^{4}}\;  \wmu_n\otimes  \wmu_n \\
S_{n43}=  \frac{1}{\sqrt n}\sum_{i=1}^n \frac{\langle X_i, \wmu_n\rangle }{\| X_i\|^{4}} \, \wmu_n \otimes X_i & \qquad &
S_{n44}=\frac{1}{\sqrt n}\sum_{i=1}^n \frac{\langle X_i, \wmu_n\rangle }{\| X_i\|^{4}}\, X_i\otimes  \wmu_n 
\\
S_{n45}=\frac{2}{\sqrt n}\sum_{i=1}^n \frac{\langle X_i, \wmu_n\rangle }{\| X_i\|^{4}}\, \langle X_i, \wmu_n\rangle\, \Gamma_i( \wmu_n)  &\qquad & S_{n46}=\, \|  \wmu_n\|^2 \, \frac{1}{\sqrt n}\sum_{i=1}^n \frac{\langle X_i, \wmu_n\rangle }{\| X_i\|^{4}}\, \Gamma_i( \wmu_n)\,.
\end{eqnarray*}
Using again the Marcinkiewicz's strong law of large numbers, we get that $n^{-2} \sum_{i=1}^n 1/{\| X_i\|^{3}}\convpp 0$, since $\esp \|X_i\|^{-3/2} <\infty$ by \textbf{A.2}. Hence,  using \textbf{A.1} and that  $\|\Gamma_i( \wmu_n)\|_\itF=1$, we get that $S_{n4j}\convprob 0$, for $j=2, 6$. On the other hand, using again that  $n^{-3/2} \sum_{i=1}^n 1/{\| X_i\|^{2}}\convpp 0$, we obtain that  $S_{n4j}\convprob 0$ for $j=3,4,5$.

It remains to study the asymptotic behaviour of $S_{n,2}$, $S_{n,3}$ and $S_{n41}$. We will show that
\begin{eqnarray}
S_{n41}- \sqrt{n}F_X(\wmu_n-\mu)&= & o_\prob(1) \label{convS41}\\
S_{n,2}+S_{n,3}- 2 \sqrt{n}S_X(\wmu_n-\mu)&= & o_\prob(1) \label{convS23} 
\end{eqnarray}
Let us begin by showing (\ref{convS41}). Denote as $W_i:\itH\to \itF$, the random objects in  $\itB$, defined as  $W_i(u)=\left(\langle X_i,u\rangle /\| X_i\|^4\right)\, X_i\otimes X_i$, for $u\in \itH$. It is easy to see that  $\| W_i\|_{\itB}\leq \| X_i\|^{-1}$ and assumption \textbf{A.2} guarantee that $\esp \left[\| X\|^{-1}\right]<\infty$, hence the strong law of large number on $\itB$ allows to conclude that
$$\frac{1}{n}\sum_{i=1}^n \frac{ \langle X_i,\cdot\rangle}{\| X_i\|^4}\, X_i\otimes X_i \convpp F_X\,,$$
where $F_X$ is defined in (\ref{Fx}), which together with \textbf{A.1} concludes the proof of  (\ref{convS41}).

To obtain (\ref{convS23}), note that  the strong law of large number on $\itH$ and the fact that $\esp \|X_i\|^{-1} <\infty$ imply that  $(1/n)\sum_{i=1}^n X_i/\| X_i\|^2\convpp \esp[X/\| X\|^2]$. Thus, if we define a sequence $\{\itT_n\}_{n\ge 1}$ of random    objects in $\itB$ as
$$\itT_n(u)= u \otimes \frac{1}{ n}\sum_{i=1}^n \frac{X_i}{\| X_i\|^{2}}+ \frac{1}{ n}\sum_{i=1}^n \frac{X_i}{\| X_i\|^{2}}  \otimes  u \quad \mbox{ for any } u\in \itH$$
we obtain that $\itT_n\convpp 2\,S_X$, where $S_X$ is defined in (\ref{Sx}). Hence, using \textbf{A.1}, we obtain  (\ref{convS41})  concluding the proof. \square

\noi \textsc{Proof of Corollary  \ref{asym}.1.} Note that from Theorem \ref{asym}.2 we get that 
$$\sqrt n \left(\wGamma^{\Signo}(\wmu_n)-\Gamma^{\Signo}(\mu)\right)=\sqrt n \left(\wGamma^{\Signo}(\mu)-\Gamma^{\Signo}(\mu)\right)+\sqrt n G_X(\wmu_n-\mu)+o_\prob(1).$$
 Now, the results follows immediately defining, for any fixed $v\in \itH$,  the  operators  $R_v:\itH\to \itF$ and $L_v:\itH\to \itF$   as $R_v(u)=u\otimes v$ and $L_v(u)=v\otimes u$ and using that $R_v^*(\Upsilon)=\Upsilon(v)$ and $L_v^*(\Upsilon)=\Upsilon^*(v)$. \square

\end{document}